\documentclass[a4paper]{article}
\usepackage{amssymb}

\newtheorem{lemma}{Lemma}[section]
\newtheorem{theorem}[lemma]{Theorem}
\newtheorem{corollary}[lemma]{Corollary}

\newtheorem{definition}[lemma]{Definition}

\newcommand{\C}{\mathbb{C}}
\newcommand{\R}{\mathbb{R}}
\begin{document}
\title{Generic metrics and the mass endomorphism on spin three-manifolds}
\author{Andreas Hermann}
\date{}
\maketitle
\begin{abstract}
  Let~$(M,g)$ be a closed Riemannian spin manifold. The constant term in the expansion of the Green function for the Dirac operator at a fixed point~$p\in M$ is called the mass endomorphism in~$p$ associated to the metric~$g$ due to an analogy to the mass in the Yamabe problem. We show that the mass endomorphism of a generic metric on a three-dimensional spin manifold is nonzero. This implies a strict inequality which can be used to avoid bubbling-off phenomena in conformal spin geometry.
\end{abstract}

\section{Introduction}
\label{intro}
Let $(M,g)$ be a closed Riemannian spin manifold of dimension~$n\geq2$. We fix a spin structure~$\sigma$ on~$M$. If~$\overline{g}$ is a metric in the conformal class~$[g]$ of~$g$, we denote the smallest positive eigenvalue of the Dirac operator~$D^{\overline{g}}$ on $(M,\overline{g},\sigma)$ by~$\lambda_1^+(\overline{g})$. Similarly let~$\lambda_1^-(\overline{g})$ be the largest negative eigenvalue of the Dirac operator. The two quantities
\begin{eqnarray*}
  \lambda_{\min}^{+}(M,[g],\sigma)&:=&\inf_{\overline{g}\in[g]}\lambda_1^{+}(\overline{g})\,\textrm{Vol}(M,\overline{g})^{1/n},\\
  \lambda_{\min}^-(M,[g],\sigma)&:=&\inf_{\overline{g}\in[g]}|\lambda_1^-(\overline{g})|\,\textrm{Vol}(M,\overline{g})^{1/n}
\end{eqnarray*}
are conformal invariants and have been treated in many articles, e.g., in \cite{hij86}, \cite{lo}, \cite{baer92}, \cite{amm03}. We denote by~$\mathbb{S}^n$ the sphere~$S^n$ with the standard metric. If we take the unique spin structure on~$\mathbb{S}^n$ we denote the above quantity by~$\lambda_{\min}^+(\mathbb{S}^n)$. It was proven in \cite{amm03a} that the inequalities~$\lambda_{\min}^+(M,[g],\sigma)\leq\lambda_{\min}^+(\mathbb{S}^n)=\frac{n}{2}\,\omega_n^{1/n}$ and~$\lambda_{\min}^-(M,[g],\sigma)\leq\lambda_{\min}^+(\mathbb{S}^n)$ hold, where~$\omega_n$ stands for the volume of the standard sphere~$\mathbb{S}^n$.
\newline One is interested in obtaining the strict inequalities
\begin{equation}
  \label{strict_ineq}
  \lambda_{\min}^+(M,[g],\sigma)<\frac{n}{2}\,\omega_n^{1/n}\quad\textrm{or}\quad
  \lambda_{\min}^-(M,[g],\sigma)<\frac{n}{2}\,\omega_n^{1/n}
\end{equation}
for several reasons. As a first application, if~$n\geq3$ and the first or the second one of these inequalities holds, then using Hijazi's inequality one can deduce that the Yamabe problem on~$(M,g)$ has a solution. This famous problem has been solved by \cite{aub} and \cite{sch}. Furthermore the inequalities imply the existence of a solution~$\varphi$ to the nonlinear partial differential equation~$D^g\varphi=\pm\lambda_{\min}^{\pm}(M,[g],\sigma)|\varphi|^{2/(n-1)}\varphi$. This equation is related via the spinorial Weierstrass representation to constant mean curvature surfaces in~$\R^3$ in the case~$n=2$ (see \cite{baer98}, \cite{fri98}, \cite{ks96}) and to constant mean curvature hypersurfaces in Calabi-Yau manifolds of real dimension 4 in the case~$n=3$ (see \cite{amm_prep}). Because of the exponent of the nonlinear term the corresponding Sobolev embedding is critical, which makes this PDE hard to solve. Also the inequalities (\ref{strict_ineq}) can be used to avoid bubbling-off phenomena in conformal spin geometry (see \cite{amm03prep}).
\newline With the aim of obtaining the inequalities (\ref{strict_ineq}) Ammann et al. considered in \cite{ahm} the case of a conformally flat manifold~$(M,g)$ with Ker$(D^g)=0$ and introduced for a fixed point~$p\in M$ the so called \emph{mass endomorphism in $p$}. It is an endomorphism of the fiber~$\Sigma^g_pM$ of the spinor bundle over $M$ and is defined as the constant term in an asymptotic expansion of the Green function of~$D^g$ in~$p$. The main result of \cite{ahm} is that the first (resp. second) of the strict inequalities (\ref{strict_ineq}) holds, if there is a point~$p\in M$ such that~$g$ is conformally flat on a neighborhood of~$p$ and such that the mass endomorphism in~$p$ has a positive (resp. negative) eigenvalue. We note that the result in \cite{ahm} was stated for a conformally flat manifold. However, the proof given there yields the strict inequality under the weaker condition that~$p$ has a conformally flat neighborhood. Thus, we are led to the question of whether one can find a point~$p$ with nonzero mass endomorphism. The answer is only known for very few manifolds. For example, the mass endomorphism on the flat torus and on the round sphere~$\mathbb{S}^n$ always vanishes, whereas on the projective spaces $\R P^{4k+3}$ one has points with nonzero mass endomorphism (see \cite{ahm}).
\newline In this article we examine the dependence of the mass endomorphism on the Riemannian metric. Our result is that in dimension 3 for a fixed point~$p$ the mass endomorphism in~$p$ is not zero for generic metrics. More precisely we obtain the following:
\begin{theorem}
  \label{generic_metric_theo}
  Let $M$ be a three-dimensional closed spin manifold with fixed spin structure~$\sigma$. Let~$\mathcal{M}_p(M)$ be the set of all Riemannian metrics~$g$ on~$M$, which are flat on a neighborhood of~$p\in M$ and satisfy Ker$(D^g)=0$. Then the subset of all Riemannian metrics with nonzero mass endomorphism in~$p$ is dense in~$\mathcal{M}_p(M)$ with respect to the~$C^1$-topology.
\end{theorem}
It was proven in \cite{m} that for a fixed spin-structure on a closed three-manifold~$M$ the generic metric satisfies Ker$(D^g)=0$. Furthermore we can approximate in the~$C^1$-topology any given metric on~$M$ by a sequence of metrics, which are conformally flat in a sufficiently small neighborhood of~$p$. Thus we obtain
\begin{corollary}
  Let $M$ be a three-dimensional closed spin manifold with fixed spin structure~$\sigma$. The set of all Riemannian metrics for which one of the strict inequalities {\rm (\ref{strict_ineq})} holds is dense in the set of all Riemannian metrics on~$M$ with respect to the~$C^1$-topology.
\end{corollary}

\section{Preliminaries}
\label{sec:2}
In this section, we review the definition of the mass endomorphism and its basic properties following \cite{ahm}. Let~$(M,g)$ be a closed spin manifold with Ker$(D^g)=0$. We denote by~$\Sigma^gM$ the spinor bundle over~$M$ and define~$\Sigma^gM\boxtimes\Sigma^gM^*$ as the bundle over~$M\times M$ whose fiber over~$(x,y)\in M\times M$ is given by Hom$(\Sigma^g_yM,\Sigma^g_xM)$. Let~$\Delta:=\{(x,x)|x\in M\}$ be the diagonal.
\begin{definition}
  A smooth section $G^g:M\times M\backslash\Delta\rightarrow\Sigma^gM\boxtimes\Sigma^gM^*$ that is locally integrable on~$M\times M$ is called the Green function for the Dirac operator~$D^g$ if in the sense of distributions~$D^g_x(G^g(x,p))=\delta_p\mathrm{id}_{\Sigma^g_pM}$, i.e. if for all~$p\in M$,~$\Psi\in\Sigma^g_pM$ and~$\varphi\in\Gamma(\Sigma^gM)$ we have
  \begin{displaymath}
    \int_{M\backslash\{p\}} \big\langle G^g(x,p)\Psi,(D^g)^*\varphi(x) \big\rangle \,dV_g(x)=\big\langle \Psi,\varphi(p) \big\rangle,
  \end{displaymath}
  where $\big\langle.,.\big\rangle$ denotes the inner product on $\Sigma^gM$.
\end{definition}
In the case $(M,g)=(\R^n,g_{eucl})$, we denote the Green function by $G^{eucl}$. One can check that
\begin{displaymath}
  G^{eucl}(x,p)\Psi=-\frac{1}{\omega_{n-1}}\frac{x-p}{|x-p|^n}\cdot\Psi,
\end{displaymath}
where $\cdot$ denotes Clifford multiplication on~$\Sigma^gM$. It is explained in \cite{ahm} that if~$(M,g)$ is flat on a neighborhood of a point~$p\in M$ one can choose a conformal chart and obtain a local trivialisation of the spinor bundle such that the Green function has the following expansion as~$x\to p$:
\begin{equation}
  \label{green_expansion}
  G^g(x,p)\Psi=G^{eucl}(x,p)\Psi+v^g(x,p)\Psi,
\end{equation}
where $\Psi\in\Sigma^g_pM$ and~$x\mapsto v^g(x,p)\Psi$ is a smooth spinor field on a neighborhood of $p$.
As in \cite{ahm} we define the mass endomorphism in~$p$.
\begin{definition}
  Let $(M,g)$ be a closed spin manifold, which is conformally flat on a neighborhood of~$p\in M$. Choose a metric~$\overline{g}\in[g]$, which is flat on a neighborhood of~$p$ and such that~$\overline{g}_p=g_p$. Let~$G^{\overline{g}}$ be the Green function. Then, we define the mass endomorphism in~$p$ as
  \begin{displaymath}
    m^g:\Sigma^g_pM\rightarrow\Sigma^g_pM,\qquad\Psi\mapsto v^{\overline{g}}(p,p)\Psi,
  \end{displaymath}
where $v^{\overline{g}}$ is the constant term in the above expansion.
\end{definition}
It is shown in \cite{ahm} that this definition does not depend on the choice of~$\overline{g}\in[g]$ and that~$m^g$ is linear and self-adjoint. There is an analogy in conformal geometry: the constant term of the Green function~$\Gamma(.,p)$ of the Yamabe operator in~$p$ can be interpreted as the mass of the asymptotically flat manifold~$(M\backslash\{p\},\Gamma(.,p)^{4/(n-2)}g)$ (see \cite{lp}). Therefore, the endomorphism is called mass endomorphism. We note that the main result of \cite{ahm} was stated for a conformally flat manifold. However, the proof given there yields the following:
\begin{theorem}
  Let $(M,g,\sigma)$ be a closed spin manifold of dimension~$n\geq2$ with Ker$(D^g)=0$. Assume that there is a point $p\in M$ which has a conformally flat neighborhood and that the mass endomorphism in~$p$ possesses a positive (resp. negative) eigenvalue. Then
  \begin{displaymath}
    \lambda_{\min}^+(M,[g],\sigma)\quad(\textrm{resp. }\lambda_{\min}^-(M,[g],\sigma))<\lambda_{\min}^+(\mathbb{S}^n)=\frac{n}{2}\,\omega_n^{1/n}.
  \end{displaymath}
\end{theorem}

\section{Proof of result}
\label{sec:3}
This section contains the proof of Theorem \ref{generic_metric_theo}.
\begin{definition}
  Let $\varphi$ be a smooth spinor field. The energy momentum tensor~$Q_{\varphi}$ for~$\varphi$ is a symmetric~$(0,2)$ tensor field on~$M$ given by
  \begin{equation}
    \label{energymomentum}
    Q_{\varphi}(X,Y):=\frac{1}{2}\,\textrm{Re }\big\langle \varphi,X\cdot\nabla^g_Y\varphi+Y\cdot\nabla^g_X\varphi\big\rangle.
  \end{equation}
\end{definition}
\begin{lemma}
  \label{dtm_lemma}
  Let $(M,g)$ be flat on a neighborhood of~$p\in M$. Let~$\Psi_0\in\Sigma^g_pM$ and define the spinor field~$G^g_{\Psi_0}$ on~$M\backslash\{p\}$ by~$G^g_{\Psi_0}(x):=G^g(x,p)\Psi_0$. Let~$K\subset M\backslash\{p\}$ be compact and let~$(g_t)_{t\in[0,1]}$ be a smooth family of metrics on~$M$ with~$g_0=g$ and supp$(g_t-g)\subset K$ for all~$t$. Define~$h:=\frac{dg_t}{dt}|_{t=0}$. Then, we have
  \begin{equation}
    \frac{d}{dt}\, \big\langle \Psi_0, m^{g_t}(\Psi_0) \big\rangle \big|_{t=0}
  =\frac{1}{2}\,\int_{M\backslash\{p\}} (h, Q_{G^g_{\Psi_0}})\,dV_g,
  \end{equation}
  where $(.,.)$ denotes the standard pointwise inner product of $(0,2)$ tensor fields.
\end{lemma}
{\bf Proof:}
Let $\varepsilon>0$ be such that~$B_{\varepsilon}(p)\subset M\backslash K$. Let~$\eta$: $M\rightarrow[0,1]$ be a smooth function with~$\eta|_{B_{\varepsilon/2}(p)}\equiv1$, supp$(\eta)\subset B_{\varepsilon}(p)$. On a neighborhood of~$p$ the Green function is given by the expansion (\ref{green_expansion}), i.e.,
\begin{displaymath}
  G^{g_t}(x,p)\Psi_0=G^{eucl}(x,p)\Psi_0+v^{g_t}(x,p)\Psi_0,
\end{displaymath}
where $\Psi_0\in\Sigma^g_pM$. For each~$x\in M$ we define a homomorphism~$w^{g_t}(x,p):\Sigma^{g_t}_pM\rightarrow\Sigma^{g_t}_xM$, such that
\begin{equation}
  \label{vdef}
  G^{g_t}(x,p)\Psi_0=\eta(x)G^{eucl}(x,p)\Psi_0+w^{g_t}(x,p)\Psi_0.
\end{equation}
In particular $m^{g_t}(\Psi_0)=w^{g_t}(p,p)\Psi_0$. Applying~$D^{g_t}$ we obtain
\begin{displaymath}
  \delta_p\Psi_0=D^{g_t}(\eta(x)G^{eucl}(x,p)\Psi_0)+D^{g_t}(w^{g_t}(x,p)\Psi_0).
\end{displaymath}
However, since $g=g_t$ in supp$(\eta)$, it follows that
\begin{equation}
  \label{Dv}
  \delta_p\Psi_0=D^g(\eta(x)G^{eucl}(x,p)\Psi_0)+D^{g_t}(w^{g_t}(x,p)\Psi_0).
\end{equation}
As a shorthand notation we introduce the spinor fields~$G^{g_t}_{\Psi_0}$,~$w^{g_t}_{\Psi_0}$,~$G^{eucl}_{\Psi_0}\in\Gamma(\Sigma^{g_t}M)$ given by
\begin{displaymath}
  G^{g_t}_{\Psi_0}(x):=G^{g_t}(x,p)\Psi_0,\quad w^{g_t}_{\Psi_0}(x):=w^{g_t}(x,p)\Psi_0,\quad G^{eucl}_{\Psi_0}(x):=G^{eucl}(x,p)\Psi_0.
\end{displaymath}
Using the definition of the Green function we obtain
\begin{eqnarray}
  \big\langle \Psi_0, w^{g_t}_{\Psi_0}(p) \big\rangle&=&
  \int_{M\backslash\{p\}} \big\langle G^{g_t}_{\Psi_0}, (D^{g_t})^*w^{g_t}_{\Psi_0} \big\rangle \,dV_{g_t}\nonumber\\
  &\stackrel{(\ref{vdef})}{=}&\int_{M\backslash\{p\}} \big\langle \eta G^{eucl}_{\Psi_0}, (D^{g_t})^*w^{g_t}_{\Psi_0} \big\rangle \,dV_{g_t}
  +\int_{M\backslash\{p\}} \big\langle w^{g_t}_{\Psi_0}, (D^{g_t})^*w^{g_t}_{\Psi_0} \big\rangle \,dV_{g_t}\nonumber\\
  \label{mpsi}
  &=&\int_{M\backslash\{p\}} \big\langle \eta G^{eucl}_{\Psi_0}, (D^g)^* w^g_{\Psi_0} \big\rangle \,dV_g
  +\int_M \big\langle D^{g_t}w^{g_t}_{\Psi_0},w^{g_t}_{\Psi_0} \big\rangle \,dV_{g_t}.
\end{eqnarray}
In the last step we have used the fact that~$g_t=g$ in supp$(\eta)$. By (\ref{Dv}) the term~$D^{g_t}w^{g_t}_{\Psi_0}$ is independent of~$t$ and therefore we get
\begin{equation}
  \label{dtDv}
  0=\frac{d}{dt}\,(D^{g_t}w^{g_t}_{\Psi_0})\big|_{t=0}=\Big( \frac{d}{dt}\,D^{g_t}\big|_{t=0} \Big) w^g_{\Psi_0}
  +D^g\Big( \frac{d}{dt}\,w^{g_t}_{\Psi_0}\big|_{t=0} \Big).
\end{equation}
Taking the derivative with respect to~$t$ in (\ref{mpsi}) we obtain
\begin{displaymath}
  \frac{d}{dt}\,\big\langle \Psi_0, m^{g_t}(\Psi_0) \big\rangle \big|_{t=0}
  =\int_M \big\langle D^g w^g_{\Psi_0},\frac{d}{dt}\,w^{g_t}_{\Psi_0}\big|_{t=0} \big\rangle \,dV_g
  +\int_M \big\langle D^g w^g_{\Psi_0},w^g_{\Psi_0} \big\rangle \frac{d}{dt}\,dV_{g_t}\big|_{t=0}.
\end{displaymath}
On the one hand, the derivative of the volume element vanishes on~$B_{\varepsilon}(p)$, on the other hand,~$D^gw^g_{\Psi_0}$ vanishes on~$M\backslash B_{\varepsilon}(p)$ by (\ref{Dv}). Hence, the second term is zero. Since~$D^g$ is self-adjoint, we conclude that
\begin{eqnarray}
  \frac{d}{dt}\,\big\langle \Psi_0, m^{g_t}(\Psi_0) \big\rangle \big|_{t=0}
  &=&\int_M \big\langle w^g_{\Psi_0},D^g\Big(\frac{d}{dt}\,w^{g_t}_{\Psi_0}\big|_{t=0}\Big)\big\rangle \,dV_g\nonumber\\
  \label{dtmpsi}
  &\stackrel{(\ref{dtDv})}{=}&-\int_M \big\langle w^g_{\Psi_0},\Big(\frac{d}{dt}\,D^{g_t}\big|_{t=0}\Big)w^g_{\Psi_0} \big\rangle \,dV_g.
\end{eqnarray}
Let $H:TM\rightarrow TM$ be given on each fiber~$T_xM$ by the unique~$g$-symmetric endomorphism~$H$ such that
\begin{displaymath}
  h_x(u,v)=g_x(H(u),v)
\end{displaymath}
for all~$u,v\in T_xM$. According to \cite{bg} the following formula holds for any spinor~$\Psi$
\begin{displaymath}
  \Big(\frac{d}{dt}\,D^{g_t}\big|_{t=0}\Big)\Psi=-\frac{1}{2}\,\sum_{i=1}^n e_i\cdot\nabla^g_{H(e_i)}\Psi
  +\frac{1}{4}\,(\textrm{div}_gh+\textrm{dTr}_gh)\cdot\Psi,
\end{displaymath}
where~$(e_i)_{i=1}^n$ is a local~$g$-orthonormal frame and div$_g$ denotes the divergence of a symmetric~$(0,2)$ tensor. By definition we have~$H(e_i)=\sum_{j=1}^n g(H(e_i),e_j)e_j=\sum_{j=1}^n h(e_i,e_j)e_j$ and since~$h$ is symmetric we obtain
\begin{displaymath}
  -\frac{1}{2}\,\sum_{i=1}^n e_i\cdot\nabla^g_{H(e_i)}\Psi
  =-\frac{1}{4}\,\sum_{i=1}^n\sum_{j=1}^n h(e_i,e_j)(e_i\cdot\nabla^g_{e_j}\Psi+e_j\cdot\nabla^g_{e_i}\Psi).
\end{displaymath}
From (\ref{dtmpsi}) it follows that
\begin{eqnarray*}
  \frac{d}{dt}\,\big\langle \Psi_0, m^{g_t}(\Psi_0)\big\rangle \big|_{t=0}
  &=&\frac{1}{4}\,\sum_{i=1}^n\sum_{j=1}^n\int_M h(e_i,e_j)\big\langle w^g_{\Psi_0},e_i\cdot\nabla^g_{e_j}w^g_{\Psi_0}+e_j\cdot\nabla^g_{e_i}w^g_{\Psi_0}\big\rangle \,dV_g\\
  &&{}-\frac{1}{4}\,\int_M \big\langle w^g_{\Psi_0},(\textrm{div}_gh+\textrm{dTr}_gh)\cdot w^g_{\Psi_0}\big\rangle \,dV_g.
\end{eqnarray*}
Since $\textrm{div}_gh+\textrm{dTr}_gh$ is a one-form, the second term is purely imaginary. Therefore, taking the real part and using the fact that~$m^{g_t}$ is self-adjoint and that supp$(h)\subset M\backslash B_{\varepsilon}(p)$ we get
\begin{displaymath}
  \frac{d}{dt}\,\big\langle \Psi_0, m^{g_t}(\Psi_0)\big\rangle\big|_{t=0}
  =\frac{1}{2}\,\int_M (h, Q_{w^g_{\Psi_0}}) \,dV_g
  =\frac{1}{2}\,\int_{M\backslash B_{\varepsilon}(p)} (h, Q_{w^g_{\Psi_0}}) \,dV_g.
\end{displaymath}
Since~$w^g_{\Psi_0}=G^g_{\Psi_0}$ on~$M\backslash B_{\varepsilon}(p)$ we finish the proof by taking the limit~$\varepsilon\to0$.
\hfill$\Box$\vspace{0.4cm}\newline
{\bf Remark:} Let~$g$ and~$\overline{g}$ be two conformally related metrics on a closed Riemannian spin manifold, i.e., there is a smooth function~$f:M\rightarrow(0,\infty)$ such that~$g=f^2\overline{g}$. According to \cite{hit74}, \cite{hij86}, \cite{bg} there is a bundle morphism~$\beta:\Sigma^gM\rightarrow\Sigma^{\overline{g}}M$ which is a fiberwise isomorphism, such that the Dirac operators on the spinor bundles~$\Sigma^gM$ and~$\Sigma^{\overline{g}}M$ are related by the formula
\begin{equation}
  \label{dirac_conform}
  D^{\overline{g}}(\beta(\psi))=f\beta(D^g(\psi))
\end{equation}
for all~$\psi\in\Gamma(\Sigma^gM)$. Moreover~$f^{-(n-1)/2}\beta$ is a pointwise isometry with respect to the inner products on~$\Sigma^gM$ and~$\Sigma^{\overline{g}}M$. Furthermore the map~$b:(TM,g)\rightarrow(TM,{\overline{g}})$ sending~$v\in T_xM$ to~$fv$ is an isometry. We also have~$\beta(X\cdot\psi)=b(X)\overline{\cdot}\beta(\psi)$ where~$\cdot$ and~$\overline{\cdot}$ denote Clifford multiplication on~$\Sigma^gM$ and~$\Sigma^{\overline{g}}M$. The connections on the spinor bundles~$\Sigma^gM$ and~$\Sigma^{\overline{g}}M$ are related by the formula
\begin{equation}
  \label{nabla_conform}
  \nabla^{\overline{g}}_X (f^{-(n-1)/2}\beta(\psi))=f^{-(n-1)/2}\beta\big( \nabla^g_X\psi+\frac{1}{2f}\,X\cdot\nabla f\cdot\psi +\frac{1}{2f}\,X(f)\psi\big)
\end{equation}
for all~$\psi\in\Gamma(\Sigma^gM)$ and all~$X\in TM$ (see \cite{lm} p. 134). We see from (\ref{dirac_conform}) that the spaces of harmonic spinors on conformally related closed Riemannian spin manifolds are isomorphic. We say that the space of harmonic spinors is conformally invariant. We also deduce that the equation~$Q_{\psi}=0$ is conformally invariant in the following sense (see also \cite{m}):
\begin{lemma}
  \label{Q_conform}
  Let $X$, $Y\in TM$,~$\psi\in\Gamma(\Sigma^gM)$. Then
  \begin{displaymath}
    Q_{\beta(\psi)}(X,Y)=f^{n-2}Q_{\psi}(X,Y).
  \end{displaymath}
\end{lemma}
{\bf Proof:} We have
\begin{eqnarray*}
  &&\big\langle Y\overline{\cdot}\nabla^{\overline{g}}_X \beta(\psi),\beta(\psi)\big\rangle\\
  &=&f^{(n-1)/2}\,\big\langle Y\overline{\cdot}\nabla^{\overline{g}}_X(f^{-(n-1)/2}\beta(\psi)),\beta(\psi)\big\rangle
  +\frac{n-1}{2}\,f^{-1}X(f)\,\big\langle Y\overline{\cdot}\beta(\psi),\beta(\psi)\big\rangle\\
  &\stackrel{(\ref{nabla_conform})}{=}&f^{(n-1)/2}\,\big\langle Y\overline{\cdot} f^{-(n-1)/2}\beta\big( \nabla^g_X\psi+\frac{1}{2f}\,X\cdot\nabla f\cdot\psi +\frac{1}{2f}\,X(f)\psi\big),\beta(\psi)\big\rangle\\
  &&{}+\frac{n-1}{2}\,f^{-1}X(f)\,\big\langle Y\overline{\cdot}\beta(\psi),\beta(\psi)\big\rangle\\
  &=&f^{-1}\,\big\langle\beta\big( Y\cdot\nabla^g_X\psi+\frac{1}{2f}\,Y\cdot X\cdot\nabla f\cdot\psi +\frac{1}{2f}\,X(f)Y\cdot\psi\big),\beta(\psi)\big\rangle\\
  &&{}+\frac{n-1}{2}\,f^{-2}X(f)\,\big\langle \beta(Y\cdot\psi),\beta(\psi)\big\rangle\\
  &=&f^{n-2}\,\big\langle Y\cdot\nabla^g_X\psi+\frac{1}{2f}\,Y\cdot X\cdot\nabla f\cdot\psi +\frac{1}{2f}\,X(f)Y\cdot\psi,\psi\big\rangle\\
  &&{}+\frac{n-1}{2}\,f^{n-3}X(f)\,\big\langle Y\cdot\psi,\psi\big\rangle\\
  &=&f^{n-2}\,\big\langle Y\cdot\nabla^g_X\psi,\psi\big\rangle
  +\frac{1}{2}\,f^{n-3}\,\big\langle Y\cdot X\cdot\nabla f\cdot\psi,\psi\big\rangle
  +\frac{n}{2}\,f^{n-3}X(f)\,\big\langle Y\cdot\psi,\psi\big\rangle.
\end{eqnarray*}
Using the fact that $Y\cdot X+X\cdot Y=-2g(X,Y)$ and that $\big\langle \alpha\cdot\psi,\psi\big\rangle$ is purely imaginary for any one-form $\alpha$, we obtain the above equation.
\hfill$\Box$\vspace{0.4cm}\newline
{\bf Proof of Theorem \ref{generic_metric_theo}:} Assume dim~$M=3$. Let~$g\in\mathcal{M}_p(M)$ with~$m^g=0$. We assume that the claim is wrong. Then there is a neighborhood~$\mathcal{U}\subset\mathcal{M}_p(M)$ of~$g$ with respect to the~$C^1$-topology such that~$m^{\tilde{g}}=0$ for all~$\tilde{g}\in\mathcal{U}$. We choose an open set~$V\subset M$ with~$p\not\in\overline{V}$ and a metric~$\tilde{g}$ which is not conformally flat on~$V$ and by choosing~$V$ sufficiently small we may assume that~$\tilde{g}\in\mathcal{U}$. In the following, we write~$g$ instead of~$\tilde{g}$. Let~$K\subset M\backslash\{p\}$ be compact and let~$(g_t)_{t\in[0,1]}$ be a smooth family of metrics on~$M$ with~$g_0=g$ and supp$(g_t-g)\subset K$ for all~$t$. Define~$h:=\frac{dg_t}{dt}|_{t=0}$. From our assumption that the claim is wrong we conclude~$m^{g_t}=0$ for all~$t$ near~$0$, hence~$\int_{M\backslash\{p\}}(h, Q_{G^g_{\Psi_0}}) dV_g=0$ for all symmetric~$(0,2)$ tensors~$h$ with supp$(h)\subset K$ by Lemma \ref{dtm_lemma}. Taking~$h=\eta Q_{G^g_{\Psi_0}}$ with suitable compactly supported functions~$\eta$ we conclude that~$Q_{G^g_{\Psi_0}}=0$ on~$K$, and since~$K$ is arbitrary we have~$Q_{G^g_{\Psi_0}}=0$ on~$M\backslash\{p\}$. As a shorthand notation we write~$\psi:=G^g_{\Psi_0}$. We define
\begin{displaymath}
  N:=\{x\in M\backslash\{p\}|\psi(x)=0\}.
\end{displaymath}
There is an open neighborhood~$B$ of~$p$ in~$M$ such that~$B\subset M\backslash N$. This can be seen from (\ref{vdef}), since~$w^g_{\Psi_0}$ is bounded on~$M$, whereas~$|G^{eucl}_{\Psi_0}(x)|_g$ becomes arbitrarily large as~$x\rightarrow p$, where~$|.|_g$ denotes the norm induced by the inner product on~$\Sigma^gM$. Hence~$W:=M\backslash(N\cup\{p\})$ is an open subset of~$M$. We replace the metric~$g|_W$ on~$W$ by a conformally equivalent metric~$\overline{g}$ with~$g=f^2\overline{g}$ for~$f=|\psi|_g^{-2/(n-1)}$. By the above remark~$f^{-(n-1)/2}\beta$ is a pointwise isometry and therefore we obtain~$|\beta(\psi)|_{\overline{g}}=1$ on~$W$. By (\ref{dirac_conform}) we have~$D^{\overline{g}}(\beta(\psi))=0$ on~$W$ and by Lemma \ref{Q_conform} we obtain~$Q_{\beta(\psi)}=0$ on~$W$. In the following, we write~$\varphi$ instead of~$\beta(\psi)$.
\newline Since dim~$W=3$ the space of spinors is~$\Sigma_3\cong\C^2$. Let~$(e_i)_{i=1}^3$ be a local~$\overline{g}$-orthonormal frame defined on an open subset~$S\subset W$. Then for every~$x\in S$ the system
\begin{displaymath}
  \mathcal{B}_x:=\{\varphi(x),e_1\overline{\cdot}\varphi(x),e_2\overline{\cdot}\varphi(x),e_3\overline{\cdot}\varphi(x)\}
\end{displaymath}
is a real basis of~$\Sigma^{\overline{g}}_xW$ and there exist functions~$a^j_i\in C^\infty(S,\R)$,~$1\leq i\leq 3$,~$0\leq j\leq3$, such that on~$S$ we have
\begin{displaymath}
  \nabla^{\overline{g}}_{e_i}\varphi=a^0_i\varphi+\sum_{j=1}^3a^j_ie_j\overline{\cdot}\varphi,\quad 1\leq i\leq3.
\end{displaymath}
Since~$0=\partial_{e_i}\big\langle\varphi,\varphi\big\rangle=2\,\mathrm{Re}\,\big\langle\varphi,\nabla^{\overline{g}}_{e_i}\varphi\big\rangle$ it follows that~$a^0_i=0$ for all~$i$. We conclude that there is a fiberwise endomorphism~$A:TW\rightarrow TW$ such that for all~$X\in TW$ we get~$\nabla^{\overline{g}}_X\varphi=A(X)\overline{\cdot}\varphi$. We will prove now that in each fiber the endomorphism~$A:T_xW\rightarrow T_xW$ is symmetric with respect to~$\overline{g}$ using an observation by Ammann. We abbreviate~$\varphi:=\varphi(x)$.
\newline If we polarize the identity~$\big\langle X\overline{\cdot}\varphi,X\overline{\cdot}\varphi\big\rangle=\overline{g}(X,X)\big\langle\varphi,\varphi\big\rangle$ we obtain
\begin{displaymath}
  \mathrm{Re}\,\big\langle X\overline{\cdot}\varphi,Y\overline{\cdot}\varphi\big\rangle=\overline{g}(X,Y)\big\langle\varphi,\varphi\big\rangle.
\end{displaymath}
Hence, for all~$j$,~$k$ with~$j\neq k$, we have~$\big\langle e_j\overline{\cdot}e_k\overline{\cdot}\varphi,\varphi\big\rangle\in i\R$. Writing~$e_1\overline{\cdot} e_2\overline{\cdot} e_3\overline{\cdot}\varphi$ as a linear combination
\begin{displaymath}
  e_1\overline{\cdot} e_2\overline{\cdot} e_3\overline{\cdot}\varphi=a_0\varphi+\sum_{j=1}^3a_je_j\overline{\cdot}\varphi
\end{displaymath}
with~$a_j\in\R$ we can take the scalar product with~$e_j\overline{\cdot}\varphi$ and considering the real part we obtain~$a_j=0$, $j\in\{1,2,3\}$. It follows that~$e_1\overline{\cdot} e_2\overline{\cdot} e_3\overline{\cdot}\varphi=a_0\varphi$ where~$a_0\in\{\pm1\}$. By possibly changing the order of the vectors~$e_1$,~$e_2$,~$e_3$ we can take~$e_1\overline{\cdot} e_2\overline{\cdot} e_3\overline{\cdot}\varphi=\varphi$.
We calculate
\begin{eqnarray*}
  \overline{g}(A(e_2),e_1)&=&\mathrm{Re}\,\langle A(e_2)\overline{\cdot}\varphi,e_1\overline{\cdot}\varphi\rangle
  =\mathrm{Re}\,\langle\nabla^{\overline{g}}_{e_2}\varphi,e_1\overline{\cdot}\varphi\rangle
  =\mathrm{Re}\,\langle e_2\overline{\cdot}\nabla^{\overline{g}}_{e_2}\varphi,\underbrace{e_2\overline{\cdot} e_1\overline{\cdot}\varphi}_{=e_3\overline{\cdot}\varphi}\rangle\\
  &=&\mathrm{Re}\,\{\langle\underbrace{D^{\overline{g}}\varphi}_{=0},e_3\overline{\cdot}\varphi\rangle
  -\langle e_1\overline{\cdot}\nabla^{\overline{g}}_{e_1}\varphi,e_3\overline{\cdot}\varphi\rangle
  -\langle e_3\overline{\cdot}\nabla^{\overline{g}}_{e_3}\varphi,e_3\overline{\cdot}\varphi\rangle\}\\
  &=&\mathrm{Re}\,\langle e_3\overline{\cdot} e_1\overline{\cdot}\nabla^{\overline{g}}_{e_1}\varphi,\varphi\rangle
  -\underbrace{\mathrm{Re}\,\langle\nabla^{\overline{g}}_{e_3}\varphi,\varphi\rangle}_{=0}
  =-\mathrm{Re}\,\langle e_2\overline{\cdot}\nabla^{\overline{g}}_{e_1}\varphi,\varphi\rangle\\
  &=&\mathrm{Re}\,\langle\nabla^{\overline{g}}_{e_1}\varphi,e_2\overline{\cdot}\varphi\rangle
  =\mathrm{Re}\,\langle A(e_1)\overline{\cdot}\varphi,e_2\overline{\cdot}\varphi\rangle=\overline{g}(e_2,A(e_1)).
\end{eqnarray*}
Similarly we obtain~$\overline{g}(A(e_1),e_3)=\overline{g}(e_1,A(e_3))$ and~$\overline{g}(A(e_2),e_3)=\overline{g}(e_2,A(e_3))$, i.e.~$A$ is symmetric with respect to~$\overline{g}$. Therefore, we can choose the basis vectors~$e_1$,~$e_2$,~$e_3$ as eigenvectors of~$A$,~$A(e_j)=\lambda_je_j$ where~$\lambda_j\in\R$. It follows that
\begin{displaymath}
  0=\mathrm{Re}\,\big\langle \varphi,e_j\overline{\cdot}\nabla^{\overline{g}}_{e_j}\varphi\big\rangle
  =\mathrm{Re}\,\big\langle \varphi,e_j\overline{\cdot} A(e_j)\overline{\cdot}\varphi\big\rangle
  =-\mathrm{Re}\,\big\langle \varphi,\lambda_j\varphi\big\rangle=-\lambda_j.
\end{displaymath}
Hence,~$\varphi$ is a parallel spinor on~$W$. By \cite{fr} the Ricci tensor on~$(W,\overline{g})$ vanishes. Since dim~$W=3$ it follows that~$(W,\overline{g})$ is flat. However, we have chosen~$g$ such that~$(M,g)$ is not conformally flat on an open set~$V$. According to \cite{baer97} the zero set~$N$ has codimension 2 at least. Therefore~$V\backslash N$ is an open set in~$W$ on which~$(W,\overline{g})$ is not conformally flat, hence not flat. This is a contradiction.
\hfill$\Box$\vspace{0.6cm}\newline
{\bf Acknowledgements: }
The author wishes to thank Emmanuel Humbert (Nancy) and Bernd Ammann (Regensburg) for many useful comments and discussions.



\end{document}